% X-sliced-and-diced-by: 'savemail' 0.1, Aug 30, 1997

\input vanilla.sty
\input definiti
\input mathchar
\input amssym.def
\input amssym.tex
\magnification 1200
\baselineskip 18pt
\def\lra{\leftrightarrow}
\def\tal{\triangleleft}
\def\e{\varepsilon}
\def\k{\kappa}
\def\a{\alpha}
\def\pmf{\par\medpagebreak\flushpar}
\def\suml{\sum\limits}
\def\lam{\lambda}
\def\lan{\langle}
\def\ran{\rangle}
\def\uhr{\upharpoonright}
\def\breta{\bar\eta}
\def\del{\delta}
\def\om{\omega}
\def\gam{\gamma}
\def\dom{\hbox{\rm dom}}
\def\pbf{\par\bigpagebreak\flushpar}

\title
More on Entangled Orders
\endtitle
\author
Ofer Shafir and Saharon Shelah\footnote"*"{This is publication no. 553 for
the second author}
\endauthor

%\centerline{This is publication number 553 for the second author}
\subheading{Introduction}

This paper grew as a continuation of [Sh462] but in the present form it can 
serve as a motivation for it as well.
We deal with the same notions, all defined in 1.1, and use just one simple 
lemma  from there whose statement and proof we repeat in 1.3.
Originally entangledness was introduced, in [BoSh] for example, in order to 
get  narrow boolean algebras and examples of the nonmultiplicativity 
of c.c-ness.
These applications became marginal when the hope to extract new such objects 
or strong colourings were not materialized, but after the pcf constructions 
which made their 
d\'ebut in [Sh-g] it seems that this notion gained independence.

Generally we aim at characterizing the existence strong and weak entangled 
orders in cardinal arithmetic terms.
In [Sh462 \S6] necessary conditions were shown for strong entangledness which 
in   a previous version was erroneously proved to be equivalent to plain 
entangledness.
 In \S1 we give a forcing counterexample to this equivalence and in \S2 we
 get those results for entangledness (certainly the most interesting case).

A new construction of an entangled order ends this section.
In \S 3 we get weaker results for positively entangledness, especially 
when supplemented with the existence of a separating point (definition 2.1).
An antipodal case is defined in 3.10 and completely characterized in 3.11.
Lastly we outline in 3.12 a forcing example showing that these two subcases of 
positive entangledness comprise no dichotomy.
The work was done during the fall of 1994 and the winter of 1995.
We thank Shani Ben David for the beautiful typing.
\subheading{\S1. Entangledness is not strong Entangledness}
\demo{Definition 1.1} (a) A linear order $(I, <)$ is called $(\mu, 
\sigma)$-entangled
if for any matrix of distinct elements from it $\langle 
t^\e_i |i<\mu, \; \; \e < \sigma_1
\rangle  (\sigma_1 < \sigma)$ and $u \subset \sigma_1$ there are 
$\a < \beta < \mu$\
satisfying $\forall \e < \sigma_1 (t^\e_\a < t^\e_\beta \leftrightarrow 
\e \in u)$.

\noindent
(b) A linear order $(I, <)$ is called $(\mu, \sigma)$-strongly entangled
 if for any matrix $\langle t^\e_i |i<\mu \; \; \e < \sigma_1 \rangle \; 
(\sigma_1 < \sigma)$ and $u \subset \sigma$, s.t. $\forall \a < \mu \forall \e_0
 \in u \forall \e_1 \in \sigma_1 \backslash u(t^{\e_0}_\a \neq t^{\e_1}_\a)$ there
 are $\a < \beta < \mu$ satisfying $\forall \e < \sigma_1 (t^\e_\a \le t^\e_\beta
\leftrightarrow \e \in u$).

\noindent
(c) A linear order is called $(\mu, \sigma)$ positively [positively$^{[*]}$] entangled if for every
 $\sigma_1 < \sigma$ and any matrix $\lan t^\e_\a | \e < \sigma_1, \a < \mu\ran$
 s.t. $\forall \e < \sigma_1 \forall \a, \beta < \mu (t^\e_\a \neq t^\e_\beta)$
 and $u \in \{ \phi, \sigma_1\}$  there are $\a < \beta [\a \neq \beta]$ satisfying 
$ \forall \e < \sigma_1 (t^\e_\a < t^\e_\beta \lra \e \in u)$.

\noindent
(d) The phrase ``$I$ is $(\mu, \sigma)$ entangled with minimal
$\mu$" stands for ``$I$ is $(\mu, \sigma)$ entangled but not $(\mu', \sigma)$
entangled for no $\mu' < \mu$''.
\proclaim{Theorem 1.2} For any cardinals $\lambda=\lambda^{<\lambda} > \theta, 
cf\mu=\k>\lambda$ there is a cardinal preserving forcing adding a $(\mu, 
\theta^+)$-entangled order  with minimal $\mu$.
In particular, it is not $(\mu, \theta^+)$ strongly entangled.
\endproclaim
\demo{Proof} Fix $\langle \mu_i|i<cf\mu\rangle$ increasing to $\mu$ and 
define $\bbp=\{p|\dom \, p\in [\mu]^{<\lambda}$, 
for some $\a < \lambda$ ran $p\subseteq^\alpha 2$, $p $ 
is $1-1$, $\forall \a (2\a \in \dom \,  p \leftrightarrow 2\a + 1 \in \dom \,  p)
\forall \a, \beta \in [\mu_i, \mu_i + \mu_i)$  $(p(2\a) < p (2\beta) \lra p
 (2\a +1) < p (2\beta + 1)) \} $
 where $<$ is the lexicographic order.
$p\le q$ iff $\dom \,  p \subset \dom \,  q$ and $\forall \a \in \dom \,  p \left(p(\a)
\tal q(\a)\right)$.
Easily $\bbp$ is $\lambda$-closed.
In order to see that it is also $\lambda^+$-c.c. (hence cardinals preserving) 
note  that  
$\forall \a < \lambda \forall p_0, p_1, p'_0, p'_1 \in^\a 2(p_0 \tal p_0'\wedge p_1
\tal p'_1\wedge p_0 \neq p_1) \to (p_0<p_1 \leftrightarrow p'_0<p'_1)$,
 so that if $\langle p_a|\a<\lambda^+\rangle$ are from $\bbp \; 
wlog \{ \dom \,  p_\a |\a < \lambda^+\}$ is a $\Delta$  system and after
``Mostowski collapsing" the function (i.e. replacing  $\{ (\a_i, \rho_i)|
i<otp \,  \dom \,  p\}$ by $\{ (i, \rho_i)| i< otp \, \dom \,  p\} )$ we get the same element.
Now define for $\a < \beta < \lambda^+$

$$
q(x) 
= \cases p_\a(x)^{\wedge} 0  &x\in \dom \,  p_\a\\
p_\beta(x)^{\wedge}  1 &x\in \dom \,  p_\beta\backslash \dom \,  p_\a
\endcases 
$$
\pmf
and this element from $\bbp$ satisfies $p_\a,  p_\beta < q$.
Any $\bbp$-generic set induces $A=\langle e_\a|\a<\mu \rangle
 \subset^\lambda 2$ which are distinct and satisfy 
$\forall \a, \beta \in [\mu_i, \mu_i +\mu_i)$ $(e_{2\a} < e_{2\beta} \lra 
e_{2 \a+1} < e_{2\beta+1})$.
Again, $A$ is ordered lexicographically.
This shows that $(A, <)$ is not $(\mu_1, 2)$-entangled for all $\mu_1<\mu$.
Suppose by contradiction that $A$ is not $(\mu, \lambda)$-entangled.
In that case there is $p\in \bbp$ and 
$p\Vdash `` \langle\underset \sim\to{t_i^\e} |i<\mu, \e<\lambda_1\rangle, 
\; \lambda_1
< \lambda,  \, u=\{2\rho|\rho<\lambda_1\}$ is a counterexample". 
For $i<\mu$ pick 
$p<p_i$ and $\langle \a(\e, i)|\e<\lambda, \rangle \subset  \dom \, p_i$ such that 
 $p_i \Vdash  `` \mathop{\bigwedge}\limits_{\e < \lambda_1} 
\underset \sim \to{t^\e_i} = e_{\a (\e,i)}"$.
Wlog the $p_i'$s have the same ``Mostowski collapse".

%$$
%q(x) = \cases p_0(x)^\cap 0 &x=\a(0, 2 \rho)\\
%p_0 (x)^\cap 1 &x=\a(0, 2\rho+1)\\ 
%p_j(x)^\cap 1 & x=\a(j, 2\rho)\\
%p_j(x)^\cap 0 &x=\a(j, 2 \rho+1)\\
%p_0(x) &\hbox{\rm else}\endcases
%
%$$

We can assume also that for any $i<\mu \{ e_{a(\e, i)} | \e < \lambda_1\}
\subset \dom\,  p_i$ and that $\langle  otp p\in \a (\e, i) | \e < \lambda_1 \rangle$
 does not depend on $i$ (here we use the inequality $k > \lambda)$.
By the $\lambda^+-c.c.$ we have two comparable elements, call them $p_i, p_j$.
Now define

$$
q(\gamma\, =\, \cases 
p_i(\gam)^\wedge 0  &\qquad \quad \gam=\a(i, 2p)\\
p_i (\gam)^\wedge 1 &\qquad \quad 
\gam =\a(i, 2p+1)\\
p_j(\gam)^\wedge 1 &\qquad \quad \gam = \a(j, 2p)\\
p_j(\gam)^\wedge 0 &\qquad \quad \gam=\a(j, 2p+1)\\
p_i(\gam)^\wedge 0 &\gam \in \dom\,  p_i \; \; \hbox{\rm not as above}\\
p_j(\gam)^\wedge 0 &\gam \in \dom\,  p_j \; \; \hbox{\rm not as above}\endcases
$$
\pbf
But $p<q \Vdash ``\underset \sim \to t$ is not a counterexample by looking at
$i, j$", a contradiction. $\square $

The following useful lemma is [Sh462 1.2(4)].
\proclaim{Lemma 1.3} If $I$ is $(\mu, 2)$ entangled then the density of $I$
is smaller than $\mu$.
\endproclaim
\demo{Proof} Otherwise define inductively a sequence of intervals $
\langle (a^0_\a, a^1_\a) | \a < \mu \rangle$  s.t.$  (a^0_\a, a^1_\a)
$
exemplifies the nondensity of
$\{ a^0_\beta, a^1_\beta|\beta<\a\}$ in $I$, i.e. disjoint to this set.
Now the matrix $\{ a^i_\a | i < 2, \a < \mu \}$ contradicts the entangledness
with respect to $u = \{ 0 \}$. $\square$

\subheading{\S2 Positive results on entangled orders}
\demo{Definition 2.1} For a linear order $I$ and $x, y \in I$ we define 
$\langle x, y \rangle \colon = (x, y)_I \cup (y, x)_I$.
We call the point $x$ $\mu$-separative if
$|\{y\in I|y<x\}|, \; |\{y\in I|y>x\} |\ge \mu$. 
Let $f(x) = \mathop{\min}\limits_{y\in I\backslash \{x\}} |\langle x, y\rangle |$.
\proclaim{Theorem 2.2}  If $\lambda = \theta^+$ is infinite and $(I, <)$ is $(\mu,
 \lambda)$-entangled linear order with minimal  such $\mu$ 
then $|\{ x \in I|f(x)<\mu\}|<\mu$.
\endproclaim
\demo{Proof} Suppose firstly $\mu$ is limit.
Fix a strictly increasing sequence of successor cardinals converging to $\mu, 
\, \langle \mu_\a|\a< cf\mu\rangle$.
Define on $I$ the equivalence relation $xEy\lra | \langle x, y \rangle |
 < \mu$.
We look for disjoint intervals $\langle I_\a | \a < cf\mu\rangle $ satisfying
$|I_\a| \ge \mu_\a$.
If any equivalence class has power smaller than $
\mu$ then choose by induction for $I_\a$ any sufficiently large but as
 yet unchosen class (there is
  such a one---otherwise all the equivalence classes are of bounded
cardinality, there are $\mu$ many and the density of $I$ is $\mu$,
contradicting 1.3).
Otherwise fix one class $J$, $|J|=
\mu$.
Pick $x \in J$ and wlog $|\{ y \in I | y > x \} | = \mu$.
Choose inductively $\langle x_\a | \a < cf\mu \rangle$ where $x_\a$ will be taken
 as any element above the previous ones and if $\a = \beta + 1$ for some $\beta $
 then $|(x_\beta, x_\a)_I|>\mu_\beta$.
Set $I_\a = (x_\a, x_{\a+1})_I$.
Next choose counterexamples for $(\mu_\a, \lambda)$-entangledness $\langle
 t^\e_i | \e < \theta \; \, i \in [\sum\limits_{\beta < \a} \mu_\beta, \mu_\a)\rangle, u = \langle 2
\e | \e<\theta\rangle$.
For any $\a < cf\mu$ choose different elements $\langle  t^\e_i | \e \in \{
 \theta, \theta+1\}$, 
 $i \in [\sum\limits_{\beta < \a} \mu_\beta, \mu_\a)\rangle
$ in $I_\a \backslash \langle t^\e_i | \e < \theta, i \in [\sum\limits_{\beta<\a}
 \mu_\beta, \mu_\a)\rangle$ (this is possible as $\theta<2^\theta \le \mu$
by the entangledness).
Wlog all the $t^\e_i$ above are with no repetitions.
This contradicts the  $(\mu, \lambda)$-entangledness with respect to $u' = 
u \cup\{ \theta \}$.
For a successor cardinal the proof is simpler since we may disregard 
the counterexamples. $\square $
\demo{Definition 2.3} For a linear order $I \;  c.c.(I)$ is the first cardinality in 
which there is no family of disjoint open intervals.
We define $h.c.c. (I) = \min \{ c.c. (J) | J \in [I]^{|I|} \}$.
\proclaim{Lemma 2.4} If $\lambda = \theta^+$ and $I$ is $(\mu, \lambda)$-entangled
 linear order with minimal such $\mu$ then for any $\{ \sigma_i | i < \theta \} 
 \subset h c.c. (I)$ we have 
$\mathop{\Pi}\limits_{i<\theta} \sigma_i < cf\mu$.
\endproclaim
\demo{Proof} Assume not. After throwing away $<\lambda$ points of $I$ we ensure that $\forall
 x \in I (f(x) = \lambda)$ (by lemma 2.2).
Suppose the theorem fails for $\{ \sigma_i | i<\theta \}$.
Choose for every $i<\theta$ a collection of disjoint intervals $\{ I^i_\a|
\a<\sigma_i\}$ and distinct functions in $\mathop{\Pi}\limits_{i<\theta} 
\sigma_i, \langle f_\a |\a<cf\mu \rangle$.
Fix counterexamples $t^\e_i, u$ and cardinals $\langle
\mu_\a|\a<cf\mu\rangle $ as above.
For any $\e<\theta$ choose $\langle t^{\theta + 2\e}_i, t^{\theta + 2\e+1}_i|i
\in [\sum\limits_{\beta<\a} \mu_\beta, \mu_\a)\rangle$
different elements from $I^\e_{f_\a (\e)} \backslash \langle t^\e_i
|\e<\theta, i \in [\sum\limits_{\beta<\a}\mu_\beta, \mu_\a)\rangle$
(remember that $|I^i_\a|=\mu)$.
Wlog all the $t^\e_i$ are with no repetitions and so $\langle t^\e_i |\e<2\theta,\; \;
 i<\mu\rangle$ contradicts the $(\mu, \lambda)$ entangledness with respect to $
u'=u\cup \langle \theta + 2\e|\e<\theta\rangle$. $\square$
\endproclaim
\proclaim{Lemma 2.5} If $\lambda=\theta^+ $ and $I$ is $(\mu, \lambda)$ entangled
linear order with minimal such $\mu$ then $\k=hc.c.(I)$ satisfies $\k^\theta \le cf\mu$.
\endproclaim
\demo{Proof} Choose $\langle \sigma_i|i<cf \k\rangle$ unbounded in $\k$..
If $cf \k\le\theta$ then $\k^\theta=\underset i<cf\kappa \to{\Pi\sigma_i} < cf\mu$ by
 lemma 2.4.
Otherwise by the same lemma $\forall \sigma < \k (\sigma^\theta < cf\mu)$ so that
 $\k^\theta = \k \cdot \underset i < cf \k\to{\sum \sigma^\theta_i} \le cf\mu$ 
(remember that $\kappa \le cf\mu)$. $\square $
\proclaim{Corollary 2.6}  If $I$ is $(\mu, \lambda) $ entangled linear 
order with  density $\chi$ then $\forall \theta < \lambda (\chi^\theta< \mu)$.
\endproclaim
\demo{Proof} Fix $\theta<\lambda$. Wlog $\lambda=\theta^+$ and $\mu$ is minimal
 for which $I$ is $(\mu, \lambda)$ entangled.
Let $\kappa=hcc(I)$.
We know that $\kappa \in \{ \chi, \chi^+\}$.
By lemma 2.5 we have to consider only the case $\kappa^\theta=\mu$.
By the proof above it follows that $cf\kappa>\theta$ and $\mu=\kappa\cdot \underset i<cf\kappa\to{
\sum\sigma^\theta_i}$ (we keep the same notation) so that
$\kappa=\mu=cf\mu$.
$\chi < \mu=\kappa$ holds by 1.3  and we can use lemma 2.4 to get
the desired conclusion.
\pmf
{\smc  Remark 2.7} The case $\chi<\kappa=\mu(=\chi^+$ follows) occurs for example in the
 construction from [BoSh$\,$210] if we assume CH.
\proclaim{Theorem 2.8} (a) If $I$ is a $(\mu, \lambda)$ entangled order with 
minimal  $\mu$ and $\lambda=\theta^+$ then 
$\lambda < h.c.c. (I) \leq cf\mu$.  \ \  
(b) If $I$ is a $(\mu, \lambda)$ entangled order with minimal $\mu$ and 
$cf\mu \neq cf\lambda < \lambda$ then $\lambda < \chi$.
\endproclaim
\demo{Proof} (a) Assume the contrary.
By 2.5 $hcc(I) \le cf\mu$ so here $\lambda \ge hcc(I)$.
For any $x \in I$ choose a strictly increasing sequence converging to it 
with minimal (hence a regular cardinal) length
$\langle a^x_\a|\a<r(x)\rangle$.
By the assumption $\forall x (r(x) <\lam)$.
As $|rang \, r|\le |\theta \cup \{ \theta\} | = \theta < \lam \le cf\mu$
 (by lemma 3.1) for some $a=\langle x_i |i<\mu \rangle \subset I$ and some $\sigma < \lam$
 \ $ \forall x \in A (r(x)=  \sigma)$.
Wlog $\forall x \in I (f(x) = \mu)$.
Define $\langle t^\e_i|i<\mu, \e < \sigma\rangle $ by induction on $i$:
for any $\e<\sigma$ choose $t^\e_i \in (a^{x_i}_\e, a^{x_i}_{\e+1})$ different from 
previously chosen $t$'s.
This contradicts the $(\mu, \lam)$ entangledness with respect to $u=\langle 2\e|\e<\sigma
\rangle$.
\
\noindent
(b)  Assume not. 
Since $cf\mu \neq cf\lambda$ for some $\theta^+<\lambda\; \, I$ is
 $(\mu, \theta^+)$ entangled with minimal $\mu$ hence we get the conclusion of
Theorem 2.2. Now in the proof of (a) $r$ is into $\lam$ since $h.c.c. (I) 
\le \chi^+$ and 
 we can ensure only its boundedness on a large $A \subset I$.
Now take $t^\e_i$ to be in $(a^{x_i}_{\e_1}, a^{x_i}_{\e_1+1})$ where $\e_1$
 is $\e$ modulo $r(x_i)$.
$\square $
\proclaim{Corollary 2.9} If $I$ is $(\mu, \lambda)$ entangled with minimal $\mu$
 and $\lam=\theta^+$ then $cf\mu>\lam$.
\endproclaim
\demo{Proof} Immediate by lemma 2.5 and theorem 2.8.(a). 
$\square $

\noindent
{\smc Remark 2.10} For inaccessible $\lam$ we have a forcing example of a $(\mu,
\lam)$ entangled order with minimal $\mu$ and $cf\mu =\lam$.
For successor $\lambda$ this is tight -- see Theorem 1.2.
\proclaim{Theorem 2.11} If $ I$ is $(\mu, \lam)$-entangled then for any $\theta<\lam$
 and for any  matrix $\langle t^\e_i|\e<\theta, \; i < \mu\rangle$ of distinct
 elements there is a sequence of disjoint intervals $\langle I_\e|\e<\theta\rangle$
 such that all the $\a$'s and $\beta$'s in the definition of entangledness can be
 chosen to satisfy $\forall \e<\theta (t^\e_\a, t^\e_\beta \in I_\e)$.
\endproclaim
\demo{Proof} Suppose the theorem fails for $I$.
Wlog $\mu$  is minimal with respect to $\lam$ and $\mu$ is singular (otherwise
$\forall \theta<\lam (\chi^\theta <\mu=cf\mu)$ 
and then there is $A \subset \mu$ of power $\mu$ and disjoint intervals
 $\langle I_\e |\e<\theta\rangle $ for them 
   $\forall \e < \theta \forall \a \in 
A(t^\e_\a \in I_\e)$ which is more than enough).
Let $\langle \mu_i|i<cf\mu\rangle $ be a strictly increasing sequence of successor 
cardinals and  $\langle t^\e_i|\e<\theta, i < \mu\rangle $  a counterexample.
As $\chi^\theta<\mu$ wlog $\mu_0>\chi^\theta$ and by induction on $i$ we 
can choose
 $\langle I^i_\e |\e<\theta, i < cf\mu \rangle $ such that $\lan I^i_\e|\e<\theta\ran$
  are disjoint for all $i<cf\mu$ and $\forall i<cf\mu \exists j(i) < cf\mu
 \forall j>j(i) \, \exists \e <\theta\, |I^i_\e \cap \{ t^\e_v|v\in [\suml_{\a<j}
\mu_\a, \mu_j)\}|< \mu_j$ hence wlog $\forall i, j<cf\mu \left(i\neq j \to \exists
\e<\theta(I^i_\e\cap I^j_\e =  \emptyset)\right)$.
As $\lan t^\e_i|i<\mu, \e<\theta \ran$ is a counterexample, for any $i< cf\mu$ there is
 $u_i<\theta$ such that $\forall \a, \beta \in [\suml_{\a<i} \mu_\a$, $\mu_i) \exists
 \e<\theta \; \; t^\e_\a<t^\e_\beta \lra \e\notin u_i$.
By a previous lemma $2^\theta<cf \mu$ so wlog the $u_i$'s are the same $u$.
Now $\lan s^\e_i|i<\mu, \e < 3 \theta \ran$ defined by $s^\e_i  =
 t^\e_{3i}, s^{\theta+\e}_i = t^\e_{3i+1},
 s^{2 \theta + \e}_i = t^\e_{3i+2} (i<\mu, \e<\theta)$ contradicts the $(\mu, 
\lam)$-entangledness with respect to $u' = u \cup [\theta, 2\theta)$. $\square $
\proclaim{Theorem 2.12}  If $I$ is $(\mu, \theta^+)$-entangled with minimal $\mu$
 then there are two $\theta^+$-closed $\mu-c.c.$ posets whose product is not $\mu-c.c.$
\endproclaim
\demo{Proof} Let $\lan x_\a|\a<\mu\ran $ be distinct elements of $I$.
Denote by $\prec$ the partial order on $E=\{(x_{2\a}, x_{2\a+1})|\a<\mu\}$
which is the produce of $<_I$ with itself.
Let $A=\{ a\in  [E]^{\le \theta}|a$ is $
\prec-$chain$\}$ and $B=\{ a\in [E]^{\le \theta}|
\neg\exists x, y \in a (x<y)\}$.
$A$ and $B$ are $\theta^+$-closed and $A\times B$ is not $\mu-c.c. $ 
since $\{ (x_{\a}, x_{2\a+1}) |\a<\mu\}$ is an 
antichain in it.
If $\lan a_\a | \a<\mu\ran \subseteq [E]^{\le \theta}$ then look at any matrix $\lan
 t^i_\e|i<\mu, \e <\theta\ran$ satisfying $\forall \a < \mu \{ (t^\a_{2\e}, t^\a_{2\e+1})|
 \e<\theta \} \supset a_\a$ and apply theorem 2.11 with respect to $u=\phi$ 
to see that it is not an  $A$-antichain and with respect to $u=\{ 2\beta|
\beta < \theta\}$ to see that it is not a $B$-antichain.
That proves the theorem. $\square $
\proclaim{Theorem 2.13}  If $\lam^{<\lam} = \lam > \beth_w$ and $2^\lam=
\lam^+$ then there is a $(\lam^+, \lam)$ entangled order,
even strongly (follows as $\lambda=\lambda^{>\lambda})$.
\endproclaim
\demo{Proof} By [Sh$\,$460 3.5] there are $\lam$ disjoint stationary subsets  of $\lambda$ 
$\lan S_\a|\a<\lam \ran$ s.t. for each $\a < \lam$  $D\ell (S_\a)$ holds.
Since $2^\lam = \lam^+$ there is a cofinal and increasing sequence of 
functions $\lan f_\a|\a<\lam^+\ran$ in $(^\lam \lam, <^*)$ where $<^*$ 
means eventual  dominance.
Fix an enumeration  of all triples $(\gamma, \bar \eta, \e)$ where $\e, 
\gamma < \lam$ and $\bar \eta = \lan \eta^\a|\a<\gamma \ran \subset {}^\e \lam$ is a
 sequence of different functions, $\lan(\gamma_\a, \bar \eta_\a, \e_\a)|\a
<\lam \ran $ (remember that $\lam^{<\lam} = \lambda)$.
Now set $A=\{f\in \,{}^{\a}\lam |\exists \beta, \delta < \lam(2\delta < 
\gamma_\beta  \wedge \a\in S_\beta \wedge \eta^{2\delta}_\beta \triangleleft f\}$ 
and define $I=\lan f_\a|a<\lam^+ \ran $ and $f<_I g$ iff $(f\uhr \a \in 
A\lra f(\a)<g(\a)$ where $\a= min\{\beta < \lam|f(\beta)\neq g(\beta)\}$.
To prove that $I$ is as required let $\gamma<\lam, u \subset \gamma$ and 
$\lan f_{\a^\beta_\nu}|\beta<\lam^+, \nu < \gamma\ran$ be as in
definition 1.1.
To simplify the notation we write $f^\beta_\nu$ for $f_{\a^\beta_\nu}$.
Wlog $\lan \a^\beta_v|\beta<\lam^+ \ran$ is increasing for all $\nu<\lam$, 
$\gamma $ is an infinite cardinal and $u=\lan 2\a|\a<\gamma\ran$.
For every $\beta<\lam^+$ there is $\e(\beta)<\lam$ s.t. $\lan 
f^\beta_\nu\uhr\e(\beta)|\nu<\gamma\ran$ are distinct so that on $B\in 
[\lam^+]^{\lam^+}$ all $\e(\beta)$ are equal to some $\e^*$
 and all $\lan f^\beta_\nu\uhr \e (\beta)|\nu < \gamma \ran$ are the same, 
to be denoted by $\bar \eta^*$.
Let $\beta $ be s.t. $(\gamma, \bar \eta^*,
 \e^*) = (\gamma_\beta, \bar \eta_\beta, \e_\beta)$
and $\lan P_\a|\a\in S_\beta\ran$ exemplify $D\ell (S_\beta)$ in an 
obviously equivalent meaning where each $P_\a $ consists of less than 
$\lam$ sequences of length $\gamma$  of functions in ${}^\a \lam$.
If some $\e_0<\lam$ and $\bar \eta = \lan\bar\eta_\a|\a<\gamma\ran 
\subset   {}^{\e_0}\lam$  satisfy that for all $i<\lam^+$ and $\delta <\lam$
there is $\zeta\in B$ s.t. $\bar \eta = \lan f^\zeta_\nu \uhr \e_0| r<\gamma 
\ran$ and  $\min \{ f^\zeta_r(\e_0)|r<\gamma\}>\delta$ then we are 
clearly done (take such $\zeta$ with respect to $(0, 0)$ then such
$\zeta'$ with respect to $(\zeta, sup \{f^\zeta_v(\e, 0)|v<\gamma \}))$.
Otherwise for every $\bar \eta $ as above there are witnesses for its 
failure, $i(\bar \eta)$ and $\delta(\bar \eta)$.
Since $\lambda^{<\lambda}=\lam$ the supremum of $i(\bar \eta)$ over all 
relevant $\breta $ is less than $\lambda'$, denote it by $i^*$.
Define $\del\colon S_\beta\to \lam$ by $\del(\a) = sup\{ \del(\breta)|\breta 
\e P_\a\}<\lam$ and 
using the cofinality of the $f_\a$'s find $\zeta\in B \backslash i^*$ for
 which  $\del <^* f_\zeta\uhr S_\beta$.
Now using $D\ell(S_\beta)$ there is $\a \in S_\beta$
 s.t. $\lan f^\zeta_v\uhr\a|v<\gamma\ran \in P_\a$, moreover we can get
$\a>sup \min \{ \e\in S\beta|\del(\e)>f^\zeta_v(\e)\}$ so $min f^\zeta_v(\a)
> \del (\a)\ge \del(\lan f^\zeta_v\uhr\a\ran)$ a contradiction. $\square $

\subheading{\S3  Results on Positively entangled orders}
\proclaim{Theorem 3.1}  If $\mu$ is minimal s.t. $I$ is $(\mu,\lam)$ [positively$^x$]
[positively] entangled then $cf\mu \ge cf\lam$.
\endproclaim
\demo{Proof} Suppose not.
We deal with positive entangledness (the other case is similar).
Fix $\lan \mu_\a|\a<cf\mu\ran$ increasing to $\mu$ and $ \lan \lam_\a|
\a<cf\mu\ran $ s.t. $I$ is not $(\mu_i, \lam^+_i)$ positively entangled and
 counterexamples $\lan t^\e_i|i<\mu_\a\; \; \e\in[\suml_{j<i}\lam_j, 
\suml_{j\le i}\lam_j)\ran$,
wlog all with respect to $u=\emptyset$.
In each row $\e$ choose fillers $\lan t^\e_i|\mu_\a \le i < \mu\ran$ 
different  from $\lan t^\e_i|i<\mu_\a\ran$.
This contradicts the $(\mu, \lam)$-positively entangledness with respect to 
$u=\phi$.
The proof for the positively$^*$ entangledness case is similar.$\square$
\proclaim{Lemma 3.2} If a $(\mu, \lambda)$ positively$^*$ entangled linear 
order $I$ has a $\mu$-separative point then $\forall \theta < \lam
(2^\theta<\mu)$.
\endproclaim
\demo{Proof} Let $x$ be such a point and suppose by contradiction 
$\theta<\lam, \; 2^\theta
\ge \mu$.
Find distinct functions $\lan f_\a|\a<\mu\ran \subset^\theta 2$.
Define $\lan t^\e_i | \e<\theta\; i < \mu\ran $ inductively in $i$ : by 
induction on 
 $\e$ choose $x_0 < x<x_1$,
different from previously chosen $t$'s  and put $t^{2\e+\ell}_i = x_\ell$ for 
$\ell \in \{0, 1\}$ if $f_i(\e)=0$ and $t^{2\e+\ell}_i = x_{1-\ell}$ else.
This contradicts the $(\mu, \lam)$-positively$^*$ entangledness. $\square$
\proclaim{Corollary 3.3} (a) If $I$ is $(\mu, \lam)$ positively$^*$ 
entangled then  $\chi=den \, I \ge \lam$.
\
\noindent
(b) If $I$ is $(\mu, \lam)$ positively$^*$ entangled then it is not 
$(\lam, 2)$ entangled.
\endproclaim
\demo{Proof} Assume $I$ is a counterexample for (a).
Wlog $\mu$ is the minimal cardinal s.t. $I$ is $(\mu, \chi^+)$ 
positively$^*$ entangled.
If there is no $\mu$-separating point in $I$ we can define inductively a 
monotone sequence in $I$ of length $cf\mu$  which is greater than $\chi$ by 
theorem 3.1, a contradiction.
If there is a $\mu$-separating point then by lemma 2.2 \ $2^\chi <\mu$,
a contradiction.
\
\noindent
(b) follows from (a) and [Sh462  1.2(4)]. $\square $
\proclaim{Theorem 3.4}  If $I$ is $(\mu, \lam)$ positively entangled then
 $\forall  \theta <\lam (2^\theta <\mu)$.
\demo{Proof} Wlog $\lam = \theta^+$.
In view  of lemma 3.2 we can assume that $I$ has no $\mu$-separating point.
It follows that $cf\mu <  \mu$.
For any $\mu_1 < \mu$ there is a $\mu_1$-separating point, 
otherwise wlog $\forall
 x \in I \; |\{ y \in I|y<x\} |<\mu$, so we can define an increasing 
sequence  of length $\mu_1+1$ and get a contradiction.
By lemma 3.2 \  $I$ is not $(\mu, \lam)$ positively$^*$ entangled for 
every $\mu_1  < \mu$.
But now if $\lan \mu_\a |\a< cf\mu \ran$ are increasing  to $\mu$, 
$\lan \lan t^\e_i|\e<\theta \; \; i\in [\mu_\a, \mu_{\a+1} \ran|\a 
< cf\mu \ran$ are counterexamples for $(\mu_\a, \lam)$ positively$^*$ 
entangledness and $\lan I_\a| \a < cf\mu\ran$ is an inductively chosen
 monotone sequence of intervals s.t. $|I_\a| \ge 
\mu_\a$ then pick for every $\a < cf\mu$ different  $\lan t^\theta_i 
| i\in [\mu_\a, \mu_{\a+1} \ran$ from $I_\a $ to contradict the 
$(\mu, \lam)$ positively entangledness  with\break  $\lan t^\e_i | \e
\le \theta, \; i < \mu \ran$.  $\square $
\endproclaim 
\proclaim{Theorem 3.5} If $I$ is $(\mu, \lam)$ positively entangled with 
minimal  such $\mu$ which has a $\mu$-separating point and $\lam=\theta^+$ 
then $\forall \theta <\lam (2^\theta <cf\mu)$.
In particular $\lam \le cf\mu$.
\endproclaim
\demo{Proof} Let $x \in I$ be $\mu$-separating and assume $\theta<\lam, 
2^\theta \ge cf\mu$.
Fix distinct $\lan f_\a|\a < cf\mu \ran \subset^\theta 2$ and choose $\lan 
t^\e_i | \e<\theta \;  \; i \in [\mu_\a, \mu_{\a+1})\ran$ counterexamples 
for $(\mu_\a, \lam)$ positively entangledness, wlog all with respect to
$u=\phi$.
For every $\e<\theta$ choose by induction on $\a \; x_0 < x < x_1$ different 
from previously chosen elements and put $t^{\theta+\e}_\a = x_\ell$ for 
$\ell \in \{0, 1\}$ if $f_\beta (\e) = 0$ and $t^{\theta + \e}_\a
= x_{1-\ell}$ else (here $\beta$  is s.t. $\a \in [\mu_\beta, 
\mu_{\beta+1})$).
 $\lan t^\e_\a |\a<\mu \; \e < \theta + \theta \ran$ contradicts the $(\mu, \lam)$
positively entangledness. $\square $
\proclaim{Corollary 3.6} If $I$ is $(\mu, \lam)$ positively entangled with minimal $\mu$
 which has a $\mu$-\break separative point and $cf\mu \neq c\del \lam$ then $\forall 
\theta < \lam (2^\theta < cf\mu)$ and $\lambda < cf\mu$ .
\endproclaim
\demo{Proof} As $cf\mu \neq cf\lam$ there is $\theta_1 < \lam $ such that $I$
 is $(\mu,  \theta^+)$ entangled with minimal $\mu$ for every $\theta_1 
\le \theta< \lam$ so we can use theorem 2.5.
 Note that the possibility $\lam =cf\mu$ is excluded by the assumption.
\pmf
{\smc Definition 3.7} A linear order $I$ is called hereditarily separative if 
every $A \in [I]^{|I|}$ has a $|I|$-separative point.
Below $(*)$ is the property $cf\lam = \aleph_0$ and $pp \lam <  2^\lam$.
The existence of a strong limit cardinal which satisfies $(*)$ is not
known to be consistent.
\proclaim{Theorem 3.8} If $I$ is hereditarily separative $(\mu, \lambda)$-positively
 entangled with minimal $\mu$, $cf\mu \neq cf\lam$ and $\lam$ is not
 inaccessible and does not satisfy $(*)$ then $\forall \theta < \lam 
(\lam^\theta< cf\mu)$.
\endproclaim
\demo{Proof} Fix $\theta < \lam$.
If for some ${\theta_1<\lam}\; \,  \; \lam \le 2^{\theta_1}$ then $\lam^\theta \le
 2^{\theta+\theta_1} < cf\mu$.
Otherwise $\lam $ is strong limit and singular as it is not inaccessible.
By [Sh410 3.7] there are under the assumptions functions $\lan f_\a | \a
 < \lam^\theta \ran \subset \; {}^\theta \lam$ satisfying $\forall \a < \beta < 
\lam^\theta \exists \e < \theta (f_\a (\e) < f_\beta (\e))$.
If the theorem fails then $\lam^\theta \ge cf\mu$.
If $A$ is an equivalence class of the equivalence relation $x Ey\lra
|\lan x, y\ran_I| < \mu$ and is of cardinality $\mu$ then pick any 
$x\in A$.
Wlog $\{ y\in A \colon y >x\} |=\mu$.
Since $I$ is hereditarily separative $\{ y \in A|y>x\}$ has 
$\mu$-separative point, call it $z$.
In particular $|(x, z)_I|=\mu$ so $x {\not \!\!E} z$, a contradiction.
We conclude that any equivalence class of $E$ is of size less than 
$\mu$ which  implies that there are at least $cf\mu$ many such classes.
By corollary 3.6 $\lambda < cf\mu$ and as $\lam$ is strong limit 
$(2^\theta)^+ <\lam$.
Choosing any $(2^\theta)^+$ distinct equivalence classes of $E$ they 
inherit the order $I$ since they are convex subsets of it so by the 
Erd\"os-Rado theorem $\theta$  from them form a monotone sequence, call it $\lan
 J_a|\a<\theta\ran $.
Replacing it by $\lan J_\a'|\a<\theta\ran$ where $J_\a' = $ convex $(J_{2\a} 
\cup J_{2\a+1})$ we ensure also  $\forall \a|J_\a' |=\mu (J_\a' $ contains an interval
 between two nonequivalent points). 
Of course, this can be done for any $\tau<\lam$ instead of $\theta$.
Starting from any such, wlog, increasing sequence $\lan J_\a |\a<cf\lam\ran$ (remember
 $cf\lam<\lam)$ we fix a strictly monotone sequence of cardinals converging
 to $\lam$, $\lan \lam_\a |\a<cf\a\ran$.
Any $J_\a$ satisfies the assumptions of the theorem so it contains by the same argument
 monotone sequence of large intervals of length $\lam_\a \;  \lan 
J^\beta_\a|\beta < \lam_\a\ran$.
If in one $J_\a$ there is no increasing sequence of length $\lam_\a$ then 
starting from decreasing intervals $\lan J_{\a'} | \a < cf\mu\ran$ inside this $J_\a$ 
we can take all the sequences decreasing.
Otherwise we take them all increasing.
Concatenating them yields a monotone sequence of intervals $\lan I_\a|\a <
\lam \ran, \forall \a ( |I_\a | = \mu)$.
Now choose $\lan \mu_\a| \a < cf\mu\ran\; \,  \lan t^\a_\e | \a < \mu \e < \theta \ran$
 as in theorem 3.5.
For all $\e < \theta $ choose by induction on $\a  $ \  $t^\a_{\theta + e} \in I_{f_\beta (\e)} \backslash
 \{ t^\gamma_{\theta + \e} | \gamma < a\}$
 where $\a \in [\mu_\beta, \mu_{\beta+1})$.
This is always possible because $\forall \a (|I_\a | = \mu)$.
Now check that $\lan t^\a_\e|\a<\mu ; \e < \theta + \theta \ran$ contradicts
 the $(\mu, \lam)$-positively entangledness. $\square $
\proclaim{Corollary 3.9} If $I$ is $(\mu, \lam)$ positively entangled hereditarily
 separative linear order with density $\chi$ and 
$cf\mu \neq cf\lambda < \lambda$ and  $\lambda $ 
does not
 satisfy $(*)$ then $\chi > \lam$.
\endproclaim
\demo{Proof} Wlog $\mu$ is minimal.
Assume that $I$ is a counterexample and deduce by corollary 3.3 (a) that
$\chi=\lam$.
Fix $A \in [I]^\lam$ dense in $I$.
for every $x \in I$ find a well  ordered sequence of elements from $A$
 converging to $x$ of minimal length $\lan a^x_\a|\a<r(x)\ran$.
 By minimality $r(x)$ is always a regular cardinal hence smaller than $\lam$.
By corollary 3.6 $\lam < cf\mu$ and by theorem 3.8  $\forall \theta < \lam (\lambda^\theta
 < cf\mu)$.
Together we get $\lam^{<\lam} < cf\mu$ so there are two distinct points in $I$ 
with the same sequences, a clear contradiction. $\square $
\demo{Definition 3.10} If $\mu$ is a singular cardinal then a linear order $I$ is
 called ``of type $s_\mu$" if it contains for some (equivalently any) sequence
 of cardinals converging to $\mu$ $\lan \mu_\a|\a<cf\mu\ran$ an isomorphic copy of
 $\bigcup\limits_{\a<cf\mu} \{ \mu_\a\} \times \mu_\a$ ordered by $(\a, \beta) < (\a_1,
 \beta_1)$ iff $\a<\a_1$ or $\a=\a_1$ and  $\beta>\beta_1$.
We say that ``$s_\mu$ is $(\mu,\lam)$ positively entangled" if some (equivalently
 any) order of type $s_\mu$ has this property.
\proclaim{Theorem 3.11}  $s_\mu$ is $(\mu, \theta^+)$-positively entangled iff
 $\theta<cf\mu$ and $(cf\mu)^\theta < \mu$.
\endproclaim
\demo{Proof} Throughout the proof fix a sequence of successor  cardinals $\lan \mu_\a|\a<cf\mu\ran$
 strictly increasing to $\mu$.
First assume $(cf\mu)^\theta < \mu$  and $\theta  < cf\mu$.
Given any  $\lan t^\a_\e|\a > \mu, \e < \theta$  as in definition 1.1(c)  then,
as $(cf\mu)^\theta <\mu$, there is $A < \mu$ of cardinality $(2^\theta)^+$ for which
 if $\a, \beta \in A$ and $\e < \theta$ then $t^\a_\e$ and $t^\beta_\e$ have the same 
 first coordinate.
Now we can find $\a, \beta \in A$ satisfying $\forall \e < \theta (t^\a_\e > t^\beta_\e)$
and $\a<\beta$.
Otherwise color $[A]^2$ with $f(\{ \a, \beta\}) = \min \{ \e<\theta| t^\a_\e < 
t^\beta_\e)\}$ (here $\a < \beta)$ and using Erd\"os-Rado get a homogeneous set of
 size $\theta$ giving rise to a decreasing sequence of ordinals of this length,
 a contradiction.
To get the other condition observe that $\bigcup_{\e<\theta} \{ \a<\mu|t^\theta_\e >
 t^\a_\e\}$ is of cardinality less than  $\mu$ as it is a union of size less than
$cf\mu$ of initial segments of $s_\mu$, which is of order type $\mu$.
For any $\a$ in its complements we have $\forall \e < \theta (t^0_\e< t^\a_\e)$.
We conclude that $s_\mu$ is $(\mu, \theta^+)$-positively entangled.

Suppose $(cf\mu)^\theta \ge \mu$, hence there are distinct $\lan f_\a|\a<\mu\ran \subset \; 
{}^\theta (cf\mu)$.
Wlog $\forall \a \ge \mu_\a (\min f_\a > \a$).
For $\e < \theta \;\;  \beta = \mu_\a + \gamma < \mu_{\a+1}$ define $t^\beta_\e = 
(f_\beta(\e), \gamma) \in s_\mu$.
Now fix any $\a < cf\mu$ and choose a partition of $\mu_{\a+2}$ to $\mu_{\a+1}$ unbounded
 sets $\lan A_\delta | \delta<\mu_{\a+1} \ran$.
For any $\e < \theta $ look at the relation on $\mu_{\a+1}\backslash \mu_\a $ 
defined by $\beta <_{\e} \gamma  \lra f_\beta (\e) < f_\gamma (\e)$.
$\prec_\e$ is a partial order with no infinite decreasing sequences so we can define a 
 rank function $g_\e$ into $\mu_{\a+2}$ satisfying $\beta \prec_\e \gamma \to g_\e (\beta)<
 g_\e (\gamma)$ by $\prec_\e$-recursion: $g_\e (\beta) = min A_\beta \backslash \sup
\{ g_\e (\gamma) | \gamma <_\e \beta\}$.
For $\beta  \in \mu_{\a+1 }\backslash \mu_\a$ set $t^\beta_{\theta+\e} = \big(\a+2, g_\e (\beta)\big)$.
By the construction the $t$'s are different in each $\mu$-row.
If $\beta<\gam<\mu$ then either $\exists \a < cf\mu (\mu_\a \le \beta 
< \gam < \mu_{\a+1})$
 in this case since the $f_\a$'s are distinct there is $\e < \theta $ 
for which $ f_\beta (\e) \neq f_\gam (\e)$; or $f_\beta (\e) < 
f_\gamma (\e)$ so $t^\beta_\e < t^\gamma_\e$ or $f_\beta (\e) > f_\gamma 
(\e)$ which implies $\beta >_\e\gam, g_\e (\beta)>g_\e (\gamma)$ and 
$t^\beta_{\theta+\e} < t^\gamma_{\theta+\e}$.
We summarize that $\forall \beta < \gam < \mu \exists \e < \theta+\theta
(t^\beta_\e < t^\gamma_\e)$ which means that $s_\mu$ is not $(\mu, \theta^+)
$-positively entangled.

Finally we show that $s_\mu$ cannot be $(\mu, (cf\mu)^+)$-positively entangled.
For this part $cf\mu$ into  $cf\mu$ mutually disjoint stationary sets $\lan A_\a|
\a < cf\mu \ran$ and enumerate their elements $A_\a = \lan a^\a_i | i < cf\mu \ran$.
Wlog $\forall \a(a^\a_0 >\a)$.
For any $\e <f\mu \; \beta = \mu_\a + \gamma < \mu_{\a +1}$  set $t^\beta_\e=
(a^\a_\e, \gamma) \in s_\mu$.
These $t$'s are different in each $\mu$-row.
Now if for some $\beta<\gamma<\mu$  $\forall \e < cf\mu$ $(t^\beta_\e < t^\gamma_\e)$
holds then necessarily there are distinct $\a, \bar \a < cf\mu$ s.t.
 $\beta \in [\mu_\a, \mu_{\a+1}), \gamma \in [\mu_{\bar \a}, \mu_{\bar \a+1})$.
The function $f=\{ (a^{\bar \a}_\e, a^\a_\e)|\e < cf\mu\}$ is a one to one
regressive function with domain $A_{\bar \a}$ which is stationary - a contradiction.
$\square $

By the above theorem one can see that lemma 3.5 does not hold generally.
\proclaim{Theorem 3.12} There is a c.c.c. forcing adding a $(\aleph_\omega, \aleph_0)$
 positively entangled linear order of density $\aleph_0$ (in particular not of 
type $s_{\aleph_\omega}$) which has no $\aleph_\omega$-separative point.
\endproclaim
\demo{Proof} Fix any $n<\omega$ and define
$\bbp=\{f$ is a function, $\dom \, f\in [n \times \aleph_\omega]^{<\om}, 
ran f \subset 2^{<\om}$, if $ \aleph_m \le \a < \beta < \aleph_{m+1}$ are in 
 $\dom \, f$ then $\exists i<n\big( (i, \a), (i, \beta) \in \dom \, f \wedge f(i, \a)
 <_{\ell x} f(i, \beta)\big)\}$.
The order is $f\le g$ iff $\dom \, f \supseteq  \dom \, g $ and $\forall x \in \dom \, g \big(g(x) \tal
 f(x)\big)$. 
If $G$ is $\bbp$ generic we define $I = \bigcup_{m<\omega} m+\{ x \in 2^\om|\forall i < \om
 \exists f \in G \exists y \in n \times [\aleph_m, \aleph_{m+1}) \big(f(y) = x\uhr i\big)\}
$ after identifying $2^\om$ with Cantor set.
The rest is almost identical to the proof of theorem 1.1.
$\square $

\centerline{References}

\item{[BoSh210]} R. Bonnet and S. Shelah, {\sl Narrow Boolean Algebras},
Annals of Pure and Applied Logic {\bf 28}(1985), pp. 1--12..
\item{[Sh410]} S. Shelah, {\sl More on Cardinal Arithmetic},
Archive for Mathematical Logic {\bf 32}(1993), pp. 399--428.

\item{[Sh460]} S. Shelah, {\sl The Generalized Continuum Hypothesis Revisited},
Israel Journal of Mathematics.

\item{[Sh462]} S. Shelah, {\sl $\sigma$-Entangled Linear Orders and Narrowness of 
Products of Boolean Algebras}, Fundamenta Mathematicae (this volume).

\item{[Sh-g]} S. Shelah, {\sl Cardinal Arithmetic}, Oxford University Press, 1994.

\bye